\documentclass[12pt]{amsart}
\pdfoutput=1 
\usepackage{bbm}
\usepackage{amsmath,amssymb,amsthm,mathabx,mathtools}
\usepackage{mathrsfs}
\usepackage[backref=page, bookmarks=false, hidelinks]{hyperref}
\hypersetup{
    pdfauthor={Angel Toledo},
    pdftitle={TTCs in derived categories of varieties with big canonical bundle},
    colorlinks=true,
    citecolor=blue,
    linkcolor=black,
}
\newcommand{\dtee} {
\otimes_{X}^{\mathbb{L}}
}
\newcommand{\dteee} { \otimes^{\mathbb{L}} }

\newcommand{\F}{\mathscr{F}}

\newcommand{\I}{\mathscr{I}}
\newcommand{\Ll}{\mathscr{L}}

\newcommand{\Ox}{\mathscr{O}}

\newcommand{\T}{\mathscr{T}}

\newcommand{\der}[1][X]{D^{b}(#1)}

\theoremstyle{plain}
\newtheorem{thm}{Theorem}[section]
\newtheorem*{thm*}{Theorem} 

\newtheorem{prop}[thm]{Proposition}
\newtheorem{lemma}[thm]{Lemma}
\newtheorem*{lemma*}{Lemma}
\newtheorem{defn}[thm]{Definition}
\newtheorem*{defn*}{Definition}

\newtheorem*{conjecture*}{Conjecture}
\newtheorem{exmp}{Example}[section]
\newtheorem{obs}[thm]{Remark}
\newtheorem{cor}[thm]{Corollary}
\newtheorem*{cor*}{Corollary}

\numberwithin{equation}{section}

\begin{document}
\title[TTCs on $\der$ for varieties of general type]{Tensor triangulated category structures in the derived category of a variety with big (anti-)canonical bundle}
%\author{Angel Toledo}
\author{Angel Toledo}
\address{ Laboratoire J. A. Dieudonné UMR CNRS-UNS 7351.
Université Côte d'Azur, Nice, France}
\email{toledo@unice.fr}
\date{}
\maketitle
\begin{abstract}
    Let $X$ be a smooth projective variety over $\mathbb{C}$ with big (anti-)canonical bundle. It is known that in this situation the Balmer spectrum of the tensor triangulated category of perfect complexes $Perf(X)$ of $X$ equipped with the derived tensor product $\dtee$ recovers the space $X$. In this work we study the possible tensor triangulated category structures one can put on $Perf(X)$. As an application we prove a monoidal version of the well-known Bondal-Orlov reconstruction theorem.
\end{abstract}
\tableofcontents
\vspace{1in}
\newpage
\section{Introduction}
In \cite{bondal2001reconstruction}, Bondal and Orlov showed that if $X$ is a smooth projective variety over $\mathbb{C}$ with ample (anti-)canonical bundle then its bounded derived category $\der$ completely recovers the space. More precisely, they showed that
\begin{thm}\cite[Theorem 2.5]{bondal2001reconstruction}\label{thm:bonorvreconstruction}
Let X be an irreducible smooth projective variety with ample (anti-)canonical bundle. If $\der\simeq D^{b}(Y)$ for some other smooth algebraic variety Y, then $X \cong Y$.
\end{thm}
This theorem came in contrast with the discovery by Mukai (\cite{mukai1987fourier}) that for an abelian variety $A$, there exists an equivalence as triangulated categories $D^{b}(A)\simeq D^{b}(\hat{A})$ between the bounded derived category of $A$ and the bounded derived category of its dual $\hat{A}$. \\
This observation sparked the study of what is now called Fourier-Mukai partners of a given variety $X$, those varieties which are triangulated equivalent to the bounded derived category of $X$. \\
Bondal and Orlov's reconstruction pointed out that a (birational) geometric condition on the variety can introduce some control on these derived equivalences and with this in mind Kawamata generalized this theorem for varieties with big (anti-)canonical bundle clarifying from a geometric point of view what is the role of this condition on the possible equivalence of derived categories. Namely he showed:
\begin{thm}\cite[Theorem 1.4]{kawamata2002d}
Let $X,Y$ be smooth projective varieties such that there is an equivalence \[\mathcal{F}:D^{b}(X)\overset{\simeq}{\longrightarrow} D^{b}(Y)\] as triangulated categories, then
\begin{enumerate}
    \item dim X = dim Y.
    \item If the canonical divisor $K_{X}$ is nef, so is $K_{Y}$ and there is an equality in the numerical Kodaira dimensions $\nu(X)$ and $\nu(Y)$. 
    \item If X is of general type, then X and Y are birational and furthermore, there is a smooth projective variety $p:Z\to X$, $q:Z\to Y$ such that $p^{\ast}K_{X} \simeq q^{\ast}K_{Y}$.
\end{enumerate}
\end{thm}
This theorem should be understood as a strong indication of a relationship between the birational geometry of a variety and its derived category. \\
On the other hand, Balmer showed in \cite{balmer2002presheaves,bondal2001reconstruction} that when equipped with the derived tensor product $\dtee$, the derived category of perfect complexes $Perf(X)$ of any coherent scheme $X$ can recover the space $X$ by what is now known as the Balmer spectrum $Spc(Perf(X),\dtee)$. The Balmer spectrum can be constructed for a general tensor triangulated category, a triangulated category equipped with a compatible monoidal structure, and produce a locally ringed space. \\
The existence of non isomorphic Fourier-Mukai partners $Y$ for a smooth variety $X$ implies using the Balmer spectrum construction that the bounded derived category $D^{b}(X)$ can be equipped with at least as many tensor triangulated category structures as non-isomorphic Fourier-Mukai partners, up to monoidal equivalence. \\
In other words, if $FM(X)$ is the set of isomorphism classes of Fourier-Mukai partners of $X$ and $TTS(X)$ is the set of equivalence classes of tensor triangulated category structures on the bounded derived category $D^{b}(X)$ there exists an injection
\begin{align*}
FM(X)&\rightarrow TTS(X)\\
Y&\mapsto (\otimes_{Y}^{\mathbb{L}}, \Ox_{Y})
\end{align*}
Where the pair $(\otimes_{Y}^{\mathbb{L}}, \Ox_{Y})$ denotes the tensor triangulated category structure given by the derived tensor product $\otimes_{Y}^{\mathbb{L}}$ with unit $\Ox_{Y}$. \\
Our main interest in this work is the study of this function, its surjectivity and the properties that one can deduce about possible tensor triangulated category structures outside of the image of this injection, all under the condition that the (anti-)canonical bundle of $X$ is big. \\
In Section~\ref{sec2} we give a brief general overview of the results we will need about general derived categories of quasi-coherent sheaves on a smooth projective variety, together with a reminder of the Balmer spectrum construction through Thomason's classification theorem. \\
In Section\ref{section3}, given a tensor triangulated category structure $(\der, \boxtimes,\mathbbm{1})$ with unit $\mathbbm{1}$ on a bounded derived category $\der$, we introduce the notion of almost spanning class with respect to a thick subcategory $I$ (Definition \ref{defn:almostspanning}) and we show (Theorem \ref{thm:almostspanning}) that if $X$ is a smooth projective variety of general type then there exists a proper tensor ideal $I_{X^{\ast}}$ of $(\der,\dtee,\Ox_{X})$ such that the set of tensor powers of $\omega_{X}$ forms an almost spanning sequence with respect to this ideal $I_{X^{\ast}}$.
This result is meant to highlight the more general behavior of almost spanning classes through the use of Thomason's classification theorem and properties of the Balmer spectrum. We see that this collection of objects can be used to prove the following theorem:
\begin{lemma}(Lemma \ref{lemma:picardlemma})
Suppose $X$ is a smooth projective variety of general type. If $\boxtimes$ is a tensor triangulated structure on $D^{b}(X)$ with unit $\Ox_{X}$, and $U$ is a $\boxtimes$-invertible object  such that $U\boxtimes I_{X^{\ast}} \subseteq I_{X^{\ast}}$. Then there is a natural equivalence between the functors induced by $U\boxtimes \_ $ and $U\dtee\_$ in $D^{b}(X)/I_{X^{\ast}}$.
\end{lemma}
When the $\dtee$-tensor ideal $I_{X^{\ast}}$ is also a $\boxtimes$-tensor ideal for a tensor triangulated category structure as described in the previous lemma, then we obtain that the Picard group of $\boxtimes$-invertible objects is a subgroup of the Picard group of $\dtee$-invertible objects (Corollary \ref{cor:bigsubideal}). This hypothesis holds true in particular when the (anti-)canonical bundle of $X$ is ample. \\
With this observation, our main corollary is the following monoidal version of the Bondal-Orlov reconstruction theorem:
\begin{cor}(Corollary \ref{cor:monoidalbondalorlov})
Let $X$ be a smooth projective variety with ample (anti-)canonical bundle, then if $\omega_{X}[n]$ is an invertible object for a tensor triangulated structure $\boxtimes$ on $\der$ with unit $\Ox_{X}$ then $\boxtimes$ and $\dtee$ coincide on objects.
\end{cor}
The results in this work were obtained as part of the author's PhD thesis at the Laboratoire J.A. Dieudonné at the Université Côte d'Azur. The author would like to thank his advisor Carlos Simpson for many discussions and to Ivo Dell'Ambrogio and Bertrand Toën for their careful and valuable comments on the thesis manuscript. The PhD thesis was partially financed by the CONACyT-Gobierno Francés 2018 doctoral scholarship. 

\section{Derived categories and the Balmer reconstruction}\label{sec2}
Through the rest of this work we will be working exclusively with smooth projective varieties over $\mathbb{C}$. We will from now on omit the mention of the base field in our exposition. Some of the material presented in this section can be found in deeper detail in \cite{huybrechts2006fourier}. \\
The goal of this section is to introduce the basic results and notions we will be using for our results. \\
Let us start by recalling that if $X$ is a smooth projective variety then there exists an equivalence as triangulated categories between the derived category $Perf(X)$ of perfect complexes on $X$ and the bounded derived category $\der$. As a consequence of this whenever we work with such a variety we will at times make no distinction between these two categories. \\
One important feature of these categories is the existence of Serre functors, let us recall:
\begin{defn}
Let $\T$ be a triangulated category an autoequivalence $S:\T\to \T$  satisfying $Hom(A,B)\cong Hom(B,S(A))^{\ast}$ for all objects $A,B\in \T$, is called a Serre functor.
\end{defn}
\begin{exmp}\label{exmp:serrefunctor}
Specifically if for example the triangulated category is a derived category of a smooth projective scheme of dimension n, we have Grothendieck-Verdier duality which implies that for every pair of objects $M,N\in\der$, $Hom(M,N)=Hom(N,M\otimes \omega_{X}[n])^{\ast}$ where $\omega_{X}$ is the canonical bundle of X. 
\end{exmp}
This notion was first defined by Kapranov and Bondal in \cite{bondal1990representable}. The following two properties of the Serre functor are essential to our work:
\begin{lemma} \label{serrecommutes}(Proposition 1.3 \cite{bondal2001reconstruction})
Let $\T$ be a triangulated category with Serre functor S, and let $\psi:\T\to \T$ be any autoequivalence, then $\psi \circ S \cong S \circ \psi$.
\end{lemma}
\begin{prop}\label{thm:uniqueserrefunctor}(Proposition 3.4 \cite{bondal1990representable})
Let $\T$ be a triangulated category and let $S$ be a Serre functor in $\T$, then it is unique up to graded isomorphism.
\end{prop}
This latter proposition implies that whenever the Serre functor exists it is part of the data of the given category. In our case of interest, as one can write this functor using the derived tensor product $\dtee$ we have now some possible control on the monoidal structure $\dtee$ directly from the category without knowledge of $X$. \\
Another crucial notion we will use is that of spanning classes, we recall the definition:
\begin{defn}\label{defn:spanningclass}
A collection of objects $\{X_{i}\}\subseteq \T$ of a triangulated category is called a spanning class if:
\begin{enumerate}
    \item If $Hom(X_{i},D[j])=0$ for all $i,j$ then $D\simeq 0$
    \item If $Hom(D[j],X_i)=0$ for all $i,j$ then $D\simeq 0$
\end{enumerate}
\end{defn}
However, whenever the Serre functor exists in the triangulated category we see that only one of the conditions is necessary and the other will be automatically satisfied by use of the Serre functor isomorphism. \\
A general way to produce spanning classes in derived categories of abelian categories is from ample sequences:
\begin{defn}
We call a collection of objects of an abelian category $\mathcal{A}$, $\{ L_{i} \}\subset \mathcal{A}$ an ample sequence if the following conditions are met:
For $i<< 0$, and all $A\in \mathcal{A}$
\begin{enumerate}
    \item $Hom(L_{i},A)\otimes_{k} L_{i}\to A$ is surjective. 
    \item $Hom(A,L_{i})=0$
    \item $Ext^{j}(L_{i},A)=0$, $j\not=0$
\end{enumerate}
\end{defn}
As the name suggest, an important example of such sequences comes from collections of tensor powers of ample line bundles. The relation between the two notions of spanning class and ample sequence was shown by Bondal and Orlov in the following result:
\begin{lemma}\label{lemma:amplespanning}
Let $\mathcal{A}$ be an abelian category of finite homological dimension and let $\{L_{i}\}$ be an ample sequence, then the collection $\{L_{i}\}$ seen as objects of $\mathcal{D}(\mathcal{A})$ form a spanning class.
\end{lemma}
The following example illustrates how we should be exploiting the existence of ample sequences. 
\begin{exmp}\label{exmp:resolvebyample}
Let $X$ be a smooth projective variety with ample canonical bundle. Then the set $\{\omega_{X}^{\otimes_{X}} \}$ forms an ample sequence and so by the previous lemma it forms a spanning class in the derived category $\der$. \\
As a consequence, we see that any complex $\F$ of coherent sheaves can be resolved by tensor powers of the canonical bundle $\omega_{X}$. In other words, there exists a sequence:
\[ 0\to \oplus_{j_{0}}(\omega_{X}^{\otimes i_{0}})\to \dots \to \oplus_{j_{k}}(\omega_{X}^{\otimes i_{k}}) \to \F\to 0\]
\end{exmp}
\begin{obs}\label{obs:resolvebyautoeq}
We remark too that in general for a triangulated category $\T$ with a spanning class $\Omega\subset \T$, if $\phi:\T\to \T$ is an autoequivalence then the set $\phi(\Omega)$ is too a spanning class. \\
In the example above, this remark implies that one can resolve any complex $\F$ by tensor powers of sheaves of the form $\omega_{X}(i)[j]$ for a fixed $i,j\in \mathbb{Z}$.
\end{obs}
\subsection{Tensor triangulated geometry}
When dealing with derived categories of coherent sheaves on a variety one can equip this category with a monoidal structure given by the derived tensor product. One can axiomatize this sort of structure in what is known as a tensor triangulated category. \\
In this subsection we recall Balmer's spectrum construction which inputs a tensor triangulated category and outputs a locally ringed space which as we will see recovers a variety whenever we work with the derived category of perfect complexes on said variety.
\begin{defn}\label{defn:ttc}
A tensor triangulated category (TTC for short) $\T$ is a triangulated category together with the following data: \\
\begin{enumerate}
    \item A closed symmetric monoidal structure given by a functor $\otimes:\T\times \T\to \T$ additive and exact ( with respect to the k-linear structure ) on both entries.
    \item The internal Hom functor $\underline{hom}:\T\times\T\to \T$ sends triangles to triangles ( up to a sign ). 
    \item Coherent natural isomorphisms for each n and m in $\mathbb{Z}$, $r:x\otimes (y[n])\cong (x\otimes y)[n]$ and $l:(x[n])\otimes y \cong (x\otimes y)[n]$ compatible with the symmetry, associative and unit coherence morphisms from the symmetric monoidal category structure. (See for example \cite[Section 2.1.1]{dell2016triangulated} for the explicit diagrams).
\end{enumerate}
\end{defn}
We will refer to a TTC by the triple $(\T,\otimes,\mathbbm{1}_{\T})$ where $\otimes$ refers to the monoidal structure and $\mathbbm{1}_{\T}$ to the unit object. Often if there is no confusion or the unit plays no role we will omit it and write $(\T,\otimes)$ instead. \\
At times when we deal with a fixed underlying triangulated category $\T$ we will write $\otimes$ or $(\otimes,\mathbbm{1})$ to refer to a tensor triangulated category structure on $\T$.
Let us remark however that the functor $\otimes$ and unit $\mathbbm{1}_{\T}$ do not completely determine a tensor triangulated category since the compatibility conditions in the symmetric monoidal category structure can in principle change while maintaining the functor $\otimes$ and unit $\mathbbm{1}$. As we will explain in the following this does not represent a problem for our purposes.  \\
We proceed with a number of definitions.
\begin{defn}\label{def:thick}
Let $\T$ be a triangulated category, and $\I\subseteq \T$ a full triangulated subcategory, we say that it is thick if it is closed under direct summands. So that if $A\oplus B \in \I$ then $A,B\in \I$.
\end{defn}
\begin{defn}
Let $(\T,\otimes)$ be a TTC. We will say that a thick subcategory $\I\subset \T$ is a $\otimes-\text{ideal}$ if for every $A\in \T$ we have $A\otimes \I\subset \I$
\end{defn}
\begin{defn}
Let $(\T,\otimes)$ be a tensor triangulated category. Let $\I$ be a $\otimes$-ideal, we will say that it is prime if for any $A,B\in \T$ with $A\otimes B\in \I$ then $A\in \I$ or $B\in \I$.
\end{defn}
As in affine algebraic geometry we can define the spectrum of a tensor triangulated category. 
\begin{defn}
Let $(\T,\otimes,\mathbbm{1})$ be a tensor triangulated category, the set of all prime $\otimes$-ideals will be denoted by $Spc(\T,\otimes,\mathbbm{1}))$ (alternatively $Spc(\T)$, $Spc(\otimes,\mathbbm{1})$ or $Spc(\otimes)$ depending on which information is clear from context).
\end{defn}
Importantly, whenever the triangulated category $\T$ is non-zero we have that $Spec(\T,\otimes)\not=\emptyset$ for any tensor triangulated category structure $\otimes$ we can put on $\T$ (see \cite[Proposition 2.3]{balmer2005spectrum}). \\
To this set we will put a topology structure.
\begin{defn}
Let $(\T,\otimes, \mathbbm{1})$ be a TTC, the support of an object $A\in \T$, denoted $supp(A)$, is the set $\{\mathfrak{p}\in Spc(\T)\mid A\not\in \mathfrak{p}\}$.
\end{defn}
\begin{lemma}\cite[Lemma 2.6]{balmer2005spectrum}
The sets of the form $\mathcal{Z}(S):=\bigcap_{A\in S} supp(A)$, for a family of objects $S\subset \T$, form a basis for a topology on $Spc(\T)$.
\end{lemma}
An important result regarding this topology is the following, which restricts the kind of spaces we should be expecting from the construction.
\begin{thm}\cite[Propositions 2.15,2.18]{balmer2005spectrum}
For any TTC $(\T,\otimes, \mathbbm{1})$, the space $Spc(\T)$ is a spectral space in the sense of Hochster, meaning it is sober and has a basis of quasi-compact open subsets.
\end{thm}
Now that the topology on $Spc(\T)$ has been chosen, the next step is to equip this space with sheaf of rings which will act as the structure sheaf. \\
To a subset $Y\subset Spc(\T)$ we can assign a thick $\otimes$-ideal denoted by $\I_{Y}$ and defined as the subcategory supported on Y, meaning $\I_{Y}:=\{A\in \T\mid supp(A)\subset Y\}$. \\
Finally, with Y as above, we denote by $\mathbbm{1}_{T_{Y}}$ the image of the unit $\mathbbm{1}$ of $\T$ under the localization functor $\pi:\T\to \T/\I_{Y}$.
\begin{defn}\label{defn:structuresheafpsc}
Let $\T$ be a nonzero TTC and we define a structure sheaf $\mathcal{O}_{Spc(\T)}$ over $Spc(\T)$ as the sheaffification of the assignment $U\mapsto End(\mathbbm{1}_{T_{Z}})$ where $Z:=Spc(\T)\backslash U$, for an open subset $U\subset Spc(\T)$.
\end{defn}
It is not hard to see the assignment $Spc(F)$ respects composition of exact monoidal functors, so if $F:\T\to \T'$ is such a functor, we get a morphism of ringed spaces since for a closed $Z= Spc(\T)\backslash U$ we have $F(\I_{Z})\subset \I_{Z'}$ where $Z'=Spc(\T')\backslash Spc(F)^{-1}(U)$ which implies there is a morphism $\mathcal{O}_{\T}\to Spc(F)_{\ast}\mathcal{O}_{\mathcal{L}}$ and so $Spc:\mathbb{TTC}\to RS$ is a functor, and under nice conditions ( for example $\T$ being rigid ) this can be shown to be a functor $Spc:\mathbb{TTC}\to LRS$.  \\
With this construction in mind we can now describe the anticipated reconstruction theorem as described by Balmer. 
\begin{thm}\cite[Corollary 5.6]{balmer2005spectrum}\label{thm:homeoperf}
    Let $X$ be a quasi-compact and quasiseparated scheme. There is a homeomorphism \[f:X\overset{\cong}{\longrightarrow} Spc(Perf(X),\dtee).\]
\end{thm}
This homeomorphism follows from Thomason's classification theorem \cite[Theorem 3.15]{thomason_1997} which establishes a correspondence between certain subsets of a quasicompact and quasi-separated scheme $X$ and $\dtee$-ideals of $Perf(X)$. The following is a general version of this classification for tensor triangulated categories as presented by Balmer in \cite[Theorem 4.10]{balmer2005spectrum}
\begin{thm}\label{thm:balmerthomasonclassification}
Let $(\T,\otimes,U)$ be a TTC. Let $\mathscr{S}$ of those subsets $Y\subset Spc(\T)$ which are unions $Y=\bigcup_{i\in I} Y_{i}$ where $Y_{i}$ are closed subsets with quasi-compact complement for all $i\in I$. Let $\mathscr{R}$ be the set of radical $\otimes$-ideals of $\T$. Then there is an order-preserving bijection $\mathscr{S}\to \mathscr{S}$ given by the assignment which sends $Y$ to the subcategory $\T_{Y}:=\{A\in \T\mid supp(A)\subset Y\}$ and with inverse sending a radical $\otimes$-ideal $\I$ to the subset $S_{\I}:=\bigcup_{A\in \I} supp(A)$.
%Let $X$ be a quasicompact and quasi-separated scheme. Denote by $\mathscr{C}$ the set of $\dtee$-ideals of the derived category $Perf(X)$ of perfect complexes on X. \\
%Denote by $\mathscr{S}$ the set of subspaces $Y\subset X$ such that $Y=\bigcup Y_{\alpha}$ where the $Y_{\alpha}$ are closed subspaces such that $X\backslash Y$ is quasi-compact. \\
%There is a bijective correspondence between $\mathscr{C}$ and $\mathscr{S}$. \\
%The bijection is given by the assignment which sends a subspace $Y\in \mathscr{S}$ to the $\dtee$-ideal whose objects are those perfect complexes $E^{\bullet}$ such that $Supph(E^{\bullet}):=\bigcup Supp(H^{i}(E^{\ast}))$ is a subset of $Y$. On the other direction it sends a $\dtee$-ideal $\I\in \mathscr{S}$ to the subspace $Y=\bigcup_{E^{\bullet}\in \I} Supph(E^{\bullet})$.
\end{thm}
Here by radical $\otimes$-ideal we mean a $\otimes$-ideal $\I$ such that whenever $A^{\otimes n}$ is in $\I$ then $A$ is in $\I$. \\
In practice every $\otimes$-ideal is automatically a radical $\otimes$-ideal and it certainly depends on the monoidal structure one can put on the triangulated category $\T$. As pointed out by Balmer in \cite{balmer2005spectrum} this condition is satisfied as soon as the tensor triangulated category is rigid, meaning that every object is dualizable. \\
When $X$ is a variety the classification theorem can be specialized to a very simple form as pointed out by Rouquier in \cite{rouquier2003categories}.
\begin{thm}\label{thm:roquierthomason}
    Let $X$ be a variety, there is a correspondence between the set of closed subsets of $X$ and $\dtee$-ideals of finite type, those ideals generated by a single object. 
\end{thm}
Using the homeomorphism from Theorem \ref{thm:homeoperf} and the construction of the structure sheaf on $Spc(Perf(X),\dtee)$ from Definition \ref{defn:structuresheafpsc} we only need the following theorem to complete the reconstruction theorem of Balmer.
\begin{thm}\cite{balmer2002presheaves}\label{thm:balmerstructuresheaf}
Let $X$ be quasi-compact and quasiseparated scheme X. There is an isomorphism $\Ox_{X}\cong \Ox_{Spc(Perf(X),\Ox_{X})}$.
\end{thm}
The following proposition should inform us how localizations behave under taking $Spc$. 
\begin{prop}\cite[Proposition 3.11]{balmer2005spectrum}\label{prop:verdiertensorfunctor}
Let $\I\subset\T$ be thick $\otimes$-ideal, then the localization functor $\pi:\T\to \T/\I$ is an exact monoidal functor and induces an homeomorphism $Spc(\T/\I)\cong \{\mathfrak{p}\in Spc(\T)\mid \I\subset \mathfrak{p}\}$.
\end{prop}
In particular when combined with the classification theorem in the form of Theorem \ref{thm:roquierthomason} we see that open subvarieties $U$ of a variety $X$ are isomorphic to $Spc(Perf(X),\dtee)/\I_{Z}$ where $Z$ is the complement of $U$ in $X$.\\
We close this section with the following remark.
\begin{obs}
So far we have been dealing with tensor triangulated categories as described in the Definition \ref{defn:ttc}, meaning we require there to be a closed symmetric monoidal category structure on $\T$. However under closer inspection one sees that nowhere in the classification theorem nor in Balmer's construction one needs the full monoidal structure. \\
In fact so far we really only need the data of a functor $\otimes:\T\times \T\to \T$ covariant and exact in each variable, together with a unit object $\mathbbm{1}$ and isomorphisms corresponding to the symmetric, associative and unit conditions. In other words, if $(\T,\otimes,\mathbbm{1})$ and $(\T,\boxtimes,\mathbbm{1}')$ are two tensor triangulated categories with underlying triangulated category $\T$ such that $\otimes \simeq \boxtimes$ for every pair of objects in $\T$, and $\mathbbm{1}\simeq \mathbbm{1}'$then the Balmer spectra $Spc(\T,\otimes,\mathbbm{1})\cong Spc(\T,\boxtimes,\mathbbm{1}')$ as locally ringed spaces. The associators, unitors and braidings of the monoidal categories have no influence in the resulting space. \\
It is this that justifies our notation $(\T,\otimes,\mathbbm{1})$ as we have mentioned before. In the following we shall keep referring to tensor triangulated categories although our results apply for slightly more general but more awkward structures. 
\end{obs}

\section{TTC's and Picard groups}\label{section3}
While the Bondal-Orlov reconstruction (Theorem \ref{thm:bonorvreconstruction}) tells us that one can directly recover a smooth projective variety $X$ with ample (anti-)canonical bundle from the derived category $\der\simeq Perf(X)$,  there are plenty of smooth projective varieties which have non-isomorphic Fourier-Mukai partners, varieties $Y$ such that $\der\simeq D^{b}(Y)$, which implies that on a given derived category $\der$ there might be many nonequivalent tensor triangulated category structures. \\
However even in the case where our variety $X$ has ample (anti-)canonical bundle as in the hypothesis of the Bondal-Orlov reconstruction theorem, it is not immediate that there is only one possible tensor triangulated category structure. It is, in principle, possible that there might be one such structure $(\der,\boxtimes,\mathbbm{1})$ such that $D^{b}(Spc(\boxtimes, \mathbbm{1}))\not\simeq \der$ and so Bondal-Orlov does not apply. \\
In some sense our motivating question is whether Balmer's reconstruction implies Bondal-Orlov. In this section we will be looking into this and related ideas by exploring the possible tensor triangulated categories one can equip on $\der$ under the slightly more general hypothesis of $X$ having a big (anti-)canonical bundle. \\
We start by mentioning the following result by Liu and Sierra from \cite{liu2013recovering} that shows in particular that there are smooth projective varieties $X$ with ample anti-canonical bundle, so under the hypothesis of Bondal-Orlov, for which the derived category $\der$ admits a tensor triangulated category structure $(\boxtimes, \mathbbm{1})$ such that $Spc(\boxtimes, \mathbbm{1})\not\cong X$. \\
Recall that there are varieties $X$ that are known to have derived categories equivalent to the derived category of representations on a quiver (possibly with relations). For example, in the presence of a full strong exceptional collection $\{E_{i}\}$ then we have that $\der$ is equivalent to $D^{b}(mod-End(\bigoplus E_{i}))$, the derived category of finitely generated modules over the algebra $End(\bigoplus E_{i})$. This latter algebra on the other hand is equivalent to the path algebra of a quiver, and so we obtain an equivalence between the derived category of $X$ and the derived category of finite dimensional representations of a quiver $Q=(P_{n},E_{ij})$. \\
The important point here is that this derived category of representations of a quiver comes with a tensor triangulated category structure induced by the tensor product of representations. To recall, let $(V_{i},p_{ik})$ and $(W_{j},q_{js})$ two such representations, then the tensor product is given entry-wise: $(V_{i},p_{ik})\otimes_{rep} (W_{j},q_{js}):= (V_{i}\otimes W_{j}, p_{ik}\otimes q_{js})$.\\
Let us denote by $(D^{b}(repQ),\otimes^{\mathbb{L}}_{rep}, \mathbbm{1}_{rep})$ the resulting tensor triangulated category structure on $D^{b}(repQ)$ by deriving this tensor product and where $\mathbbm{1}_{rep}:=(k_{i},Id_{ij})$ is the representation given by putting $k$ on every vertex  and the identity morphism in each edge of the quiver. \\
Liu and Sierra consider quivers with relations satisfying a compatibility condition with the tensor product (\cite[Definition 1.2.5]{liu2013recovering}) and say that in this case the quiver has tensor relations. \\
\begin{thm}\cite[Theorem 2.1.5.1]{liu2013recovering}\label{thm:sierraspace}
Let $Q$ be a finite ordered quiver with tensor relations. Then $Spc((D^{b}(repQ),\otimes^{\mathbb{L}}_{rep}))$ is the discrete space $\{P_{n}\}$.
\end{thm}
They also describe completely the structure sheaf in this case.
\begin{thm}\cite[Theorem 2.2.4.1]{liu2013recovering}\label{thm:sierrastructuresheaf}
Let $Q$ be a finite ordered quiver with tensor relations. Then $\Ox_{Q}:=\Ox_{Spc(\otimes^{\mathbb{L}}_{rep})}$ is the constant sheaf of algebras $k$. So that for any open $W\subset Spc(\otimes_{rep}^{\mathbb{L}})$ we have $\Ox_{Q}(W)=k^{\oplus W}$.
\end{thm}
In particular, for $X=\mathbb{P}^{n}$ we have by a well-known result of Beilinson (\cite{beilinson1978coherent}) that $\der$ is equivalent to the category of representations of a quiver with $n+1$ vertices. Thus the derived category $\der$ has a tensor triangulated category structure $(\otimes^{\mathbb{L}}_{rep},\mathbbm{1}_{rep})$ such that $Spc((\otimes^{\mathbb{L}}_{rep},\mathbbm{1}_{rep}))\not\cong X=\mathbb{P}^{n}$. As $\mathbb{P}^{n}$ is a smooth projective variety with ample anti-canonical bundle, this previous result implies that the study of tensor triangulated category structures on $\der$ is not trivial even in the cases falling under the hypothesis of the Bondal-Orlov reconstruction theorem and might shed some light in the internal structure of the derived category in itself. \\ 
\\
In general the behavior of the dynamics of the Balmer spectrum and taking derived categories can be complex. As we know that the Balmer spectrum is a locally ringed space it has an abelian category of sheaves of modules which admits a tensor product then we can derive this category as usual, however the category of sheaves of modules is in general much more complicated than a category of coherent or even quasi-coherent sheaves. \\ 
\\
Having said that, let us put ourselves in the slightly more general situation of derived categories of varieties of general type. Recall a variety is of general type if its canonical bundle is big. In particular varieties with ample canonical bundle are of general type. \\
One alternative characterization of bigness for a variety is the following:
\begin{thm}\cite[Example 2.2.9]{lazarsfeld2017positivity}
A smooth projective variety is of general type if and only if, for any sheaf $\F\in Coh(X)$, there exists an integer $i_{0}$ depending on $\F$ such that the sheaf $\F\otimes_{X}\omega_{X}^{i}$ is generically globally generated for $i>> i_{0}$.
\end{thm}
As a consequence of the Kodaira lemma (cf. \cite[Prop 2.2.6]{lazarsfeld2017positivity}) we have the corollary: \\
\begin{cor}\label{cor:opengentype}
Let $X$ be a smooth projective variety of general type, then there exists an open sub-variety $X^{\ast}$ such that for any $\F\in Coh(X)$, there exists a positive integer $i_{0}$ for which for any $i>>i_{0}$, the sheaf $\F\mid_{X^{\ast}}\otimes_{X}\omega_{X^{\ast}}^{i}$ on $X^{\ast}$ is globally generated.
\end{cor}
Let us explain the previous corollary and the nature of the open sub-variety $X^{\ast}$. We recall some basic definitions.
\begin{defn}\label{defn:augmentedlinebundle}
Let $X$ be a projective variety and $\Ll$ a line bundle on $X$, the augmented base locus is the Zariski closed set
\[ B_{+}(\Ll):=\bigcap_{m\in\mathbb{N}} B(m\Ll-A) \]
Where $A$ is any ample line bundle, and for any line bundle $\Ll'$ the set $B(\Ll')$ is defined as the intersection of the base loci of multiples of the line bundle, that is
\[ B(\Ll'):= \bigcap_{m\in \mathbb{N}} Bs(m\Ll') \]
\end{defn}
In \cite{boucksom2014augmented} the following theorem characterizing the complement of the augmented base focus is proven:
\begin{thm}\label{thm:restrictionisample}
Let $\Ll$ be a big line bundle on a normal projective variety $X$ over an algebraically closed field. Then the complement $X\backslash B_{+}(\Ll)$ of the augmented base locus is the largest Zariski open subset $U\subseteq X\backslash B(\Ll)$ such that for all large and divisible $m(\Ll)\in \mathbb{Z}$ the restriction of the morphism 
\[ \phi_{m}:X\backslash B(\Ll) \dashrightarrow \mathbb{P}H^{0}(X,m\Ll) \]
to U is an isomorphism onto its image.
\end{thm}
The following couple important observations follow immediately from the definition, the fact that the augmented base locus is independent of the choice of ample line bundle, and Kodaira's decomposition of big line bundles. 
\begin{obs}\label{properideal}
\begin{enumerate}
    \item $B_{+}(\Ll)=\emptyset$ if and only if $\Ll$ is ample.
    \item $B_{+}(\Ll)\not=X$ if and only if $\Ll$ is big. 
\end{enumerate}
\end{obs}
From the remarks above and using Thomason's classification Theorem (\ref{thm:balmerthomasonclassification}) we know that as there exists a correspondence between closed subsets of the Balmer spectrum and radical tensor ideals in the tensor triangulated category then there exists a radical tensor ideal corresponding to the augmented base locus $B_{+}(\Ll)$ for any given line bundle $\Ll$. In particular the open sub-variety $X^{\ast}$ from Corollary \ref{cor:opengentype} is the complement of the augmented base locus, $X\backslash B_{+}(\omega_{X})$ 
 and corresponds to a $\dtee$-ideal generated by a single object ( using Theorem \ref{thm:roquierthomason} ) whose homological support gives back the closed subset $B_{+}(\omega_{X})$. \\
\begin{obs}
Let us denote by $I_{X^{\ast}}$ the $\dtee$-ideal corresponding to the open subvariety $X^{\ast}$. By our previous Remark \ref{properideal} we have that this ideal must be a proper $\dtee$-ideal of $\der$ and is the ideal $0$ precisely when the (anti-)canonical bundle is ample.
\end{obs}
We would like to understand the effect of the positivity of the canonical bundle ( in this case the fact that the variety is of general type ) on the tensor triangulated structure of the category. We know from Theorem \ref{thm:uniqueserrefunctor} that the Serre functor in a triangulated category is unique up to degree whenever it exists and so it is a property of the category and not extra data. In our concrete case we know furthermore that the Serre functor is isomorphic to $\_\dtee\omega_{X}[n]$ where $n\in\mathbb{N}$ is the dimension of the variety and $\omega_{X}$ is the dualizing sheaf of $X$. \\
Let us start with a definition mimicking that of spanning class:
\begin{defn}\label{defn:almostspanning}
Let $(\T,\otimes)$ be a tensor triangulated category, let $\I\subseteq \T$ be a thick subcategory and let us denote by $\pi:\T\to \T/\I$ the localization functor. We say that a collection of objects $\Omega\subset \T$ is an almost spanning class with respect to $\I$ if the following two conditions hold.
\begin{enumerate}
    \item If $X\in \T/\I$ is such that $Hom_{\T/\I}(\pi(B),X[j])=0$ for all $B\in \Omega$ and $j\in \mathbb{Z}$, then $X\cong 0$.
    \item If $X\in \T/\I$ is such that $Hom_{\T/I}(X[j],\pi(B))=0$ for all $B\in \Omega$ and $j\in \mathbb{Z}$, then $X\cong 0$.
\end{enumerate}
\end{defn}
It is immediate to see that the previous definition is equivalent to asking that the collection $\Omega$ maps through $\pi$ to a spanning class on the quotient $\T/\I$. When the thick subcategory in question is the 0 subcategory then the definition reduces to that of a spanning class as in Definition \ref{defn:spanningclass}. \\
Additionally when the triangulated category $\T/\I$ has a Serre functor, only one of the conditions in the definition is necessary as the Serre duality implies the other automatically. \\
We would like to generalize Theorem \ref{lemma:amplespanning} but for a big canonical bundle instead of an ample one and see that a big bundle induces an almost spanning class in the derived category with respect to a $\dtee$-ideal $\I$. \\
\begin{thm}\label{thm:almostspanning}
Let $X$ be a smooth projective variety of general type. Then the collection of tensor powers $(\omega^{\otimes i}_{X})_{i\in \mathbb{Z}}$ forms an almost spanning class with respect to the tensor ideal $I_{X^{\ast}}$ in the tensor triangulated category $(\der, \dtee)$.
\begin{proof}
We need to show that $\pi(\{\omega_{X}^{\otimes i}\})$ forms a spanning class in the quotient $\der/I_{X^{\ast}}$. As $I_{X^{\ast}}$ is the ideal corresponding to the open smooth subvariety $X^{\ast}$ from Corollary \ref{cor:opengentype} then we know that there is an isomorphism $Spc(\der/I_{X^{\ast}})\cong X^{\ast}$. Since $\omega_{X}$ restricted to $X^{\ast}$ is ample by the characterization of Theorem \ref{thm:restrictionisample}, we get that $\{\omega_{X}^{\dteee i}\mid_{X^{\ast}}\}$ forms an spanning class by Lemma \ref{lemma:amplespanning} of the derived category of $X^{\ast}$ which coincides with the quotient category $\der/I_{X^{\ast}}$ by Proposition \ref{prop:verdiertensorfunctor}.
\end{proof}
\end{thm}
The main key in our arguments is the fact that one can construct, as in the ample case, a resolution for any complex of coherent sheaves on $X^{\ast}$ in terms of tensor powers of the canonical bundle $\omega_{X^{\ast}}$ of $X^{\ast}$. With the advantage that one is able to have a concrete description of the derived category of this space in terms of a quotient of the derived category of the larger variety $X$.  \\
Explicitly for any complex $A$ of coherent sheaves over $X^{\ast}$ there is a resolution:
\[ \dots\to \oplus_{j_{0}}(\omega_{X^{\ast}}^{\otimes i_{0}})\to \dots \to \oplus_{j_{k}}(\omega_{X^{\ast}}^{\otimes i_{k}})\to A\to 0.\]
\\
Another thing to notice is that in the example given above for the non equivalent tensor triangulated category structures on $D^{b}(\mathbb{P}^{n})$, one immediate issue with the two given such structures was that the units were non-isomorphic. For this reason we should proceed to work with tensor triangulated categories with a fixed unit isomorphic to $\Ox_{X}$. \\
\begin{defn}\label{defn:invertible}
Let $(\T,\otimes,\mathbbm{1})$ be a TTC, an object $X\in \T$ is $\otimes$-invertible if there exists $X^{-1}\in \T$ such that $X\otimes X^{-1}\cong \mathbbm{1}$. We will denote by $Pic(\der, \boxtimes)$ the group of isomorphism classes of $\boxtimes$-invertible objects.
\end{defn}
We will make use of the following lemma:
\begin{lemma}\label{lemma:picardlemma}
Suppose $X$ is a smooth projective variety of general type of dimension $n$. If $\boxtimes$ is a tensor triangulated structure on $D^{b}(X)$ with unit $\Ox_{X}$, and $U$ is a $\boxtimes$-invertible object  such that $U\boxtimes I_{X^{\ast}} \subseteq I_{X^{\ast}}$. Then there is a natural equivalence between the functors induced by $U\boxtimes \_ $ and $U\dtee\_$ in $D^{b}(X)/I_{X^{\ast}}$.
\begin{proof}
By our previous discussion we know that any complex can be resolved in $\der/I_{X^{\ast}}$ by a resolution
\[ \dots\to \oplus_{j_{0}}(\omega_{X^{\ast}}^{\otimes i_{0}})\to \dots \to \oplus_{j_{k}}(\omega_{X^{\ast}}^{\otimes i_{k}})\to A\to 0.\]
As the Serre functor in $D^{b}(X^{\ast})$ is given by $\_\dtee\omega_{X^{\ast}}[n']$, where $n'$ is the dimension of $X^{\ast}$ and we know any exact equivalence must commute with it, if we let $U\widehat{\boxtimes}$ and $U\widehat{\dtee}$ denote the autoequivalences of $\der/I_{X^{\ast}}$ induced respectively by $U\boxtimes$ and $U\dteee$, then we have that 
\[(U\widehat{\boxtimes} \hat{A})\widehat{\dtee} \omega_{X^{\ast}}[n'] \cong U\widehat{\boxtimes} (\widehat{A}\widehat{\dtee} \omega_{X^{\ast}}[n']).\]
As $\Ox_{X}$ is a unit for both $\otimes_{X}$ and $\boxtimes$, and after shifting by $[-n']$ we deduce
\[ U\widehat{\dtee} \omega_{X^{\ast}} \cong U\widehat{\boxtimes} \omega_{X^{\ast}}. \]
From this, the exactness of $\dteee$ and $\boxtimes$, and the resolutions in terms of $\omega_{X^{\ast}}^{i}$, we obtain the isomorphisms
\[ U\widehat{\dteee}A \cong U\widehat{\boxtimes}A. \]
\end{proof}
\end{lemma}
\begin{obs}
Let us point out the slight abuse of notation of the autoequivalence $U\widehat{\dteee}$. This functor would formally be denoted by $\widehat{U}\otimes^{\mathbb{L}}_{\der/I_{X^{\ast}}}$  as it is induced by the object $\widehat{U}$ in the tensor triangulated category $(\der/I_{X^{\ast}},\otimes^{\mathbb{L}}_{\der/I_{X^{\ast}}})$, but as the only tensor ideal we are taking a quotient by in this section is $I_{X^{\ast}}$, we believe our notation is lighter without losing sight of which functors they represent. 
\end{obs}
We have the following corollary:
\begin{cor}\label{cor:subpicard}
Let $X$ be a variety of general type and let $\boxtimes$ a tensor triangulated category structure on $D^{b}(X)$ with unit $\Ox_{X}$. Then for any $\boxtimes$-invertible object $U$ such that $U\boxtimes I_{X^{\ast}}\subseteq I_{X^{\ast}}$, the equivalence $U\widehat{\boxtimes}:\der/I_{X^{\ast}}\to \der/I_{X^{\ast}}$ induced by $U\boxtimes\_$ is equivalent to an equivalence given by objects in the group $Pic(\der/I_{X^{\ast}},\widehat{\dteee})$ of invertible $\widehat{\dteee}$-objects.
\begin{proof}
From Lemma \ref{lemma:picardlemma} we have that if $U^{-1}$ is such that $U\boxtimes U^{-1} \cong \Ox_{X}$ then in the quotient $\der/I_{X^{\ast}}$, 
\[ U\widehat{\dteee} \widehat{U^{-1}}\cong U\widehat{\boxtimes}\widehat{U^{-1}} \cong \Ox_{X^{\ast}}. \]
As $(\der/I_{X^{\ast}},\widehat{\dteee})$ is a tensor triangulated category, we have that $\widehat{U}\in \der/I_{X^{\ast}}$ is a $\widehat{\dteee}$-invertible objects.
\end{proof}
\end{cor}
In Lemma \ref{lemma:picardlemma} and Corollary \ref{cor:subpicard} above, the ideal $I_{X^{\ast}}$ might not be a $\boxtimes$-tensor ideal and thus the quotient $\der/I_{X^{\ast}}$ does not necessarily carry a tensor triangulated category structure induced by $\boxtimes$. However, our result guarantees that after passing to the quotient, the equivalences induced by the functors $U\boxtimes \_$ are equivalent to equivalences given by invertible objects in $(\der/I_{X^{\ast}},\dtee)$ induced by the same object, under the condition that $I_{X^{\ast}}$ is stable by $U\boxtimes$. \\
In particular, we have:
\begin{cor}\label{cor:bigsubideal}
Let $X$ be a variety of general type and $\boxtimes$ a tensor triangulated structure on $\der$ with unit $\Ox_{X}$. If $I_{X^{\ast}}$ is a $\boxtimes$-ideal then the Picard group $Pic(\der/I_{X^{\ast}},\widehat{\boxtimes})$ is a subgroup of the Picard group $Pic(\der/I_{X^{\ast}},\widehat{\dtee})$.
\begin{proof}
The proof is as in the previous two, if $U$ is in  $Pic(\der/I_{X^{\ast}},\widehat{\boxtimes})$ then it induces an autoequivalence of $\der/I_{X^{\ast}}$ and so it commutes with the Serre functor on $D^{b}(X^{\ast})\simeq \der/I_{X^{\ast}}$. By writing a resolution for any complex $A$ in terms of direct sums of derived tensor powers of $\omega_{X^{\ast}}$ we can use the same argument than in the proof of Lemma \ref{lemma:picardlemma} and we arrive at the isomorphisms
\[ U\widehat{\dteee} A \cong U\widehat{\boxtimes} A.\]
\end{proof}
\end{cor}
\begin{obs}
Let us point out that in the results above we have chosen to work with varieties of general type, but the same argument applies to varieties with big anti-canonical bundle.
\end{obs}
The case when our variety has an ample (anti-)canonical bundle allows us to relate the Picard group of the full derived category to that of any other tensor triangulated category structure on it. \\
The following result follows from the previous argument.
\begin{cor}\label{cor:amplepicard}
Let X be a variety with ample (anti-)canonical bundle. Then if $\boxtimes$ is a tensor triangulated category structure on $D^{b}(X)$ with unit $\Ox_{X}$, the Picard group $Pic(D^{b}(X), \boxtimes)$ is isomorphic to a subgroup of $Pic(D^{b}(X),\otimes_{X})$.
\begin{proof}
We just need to notice that in this case the $\otimes_{X}$-ideal from Corollary \ref{cor:opengentype} is the 0 ideal and thus we can resolve any object $A\in D^{b}(X)$ by a sequence of powers of the Serre functor. By the same reasoning as above we see that 
\[ U\dteee A \cong U\boxtimes A.\]
\end{proof}
\end{cor}
One thing to note here is that although Bondal and Orlov had already classified the group of autoequivalences of a derived category of a variety with ample (anti-)canonical bundle, we are working without the condition of an equivalence between the derived category of the Balmer spectrum of $\boxtimes$ and the derived category $\der$, and as such it is not immediate from their result that the Picard group of $\boxtimes$ must involve invertible sheaves over $X$. \\
In other words, as $Spc(\boxtimes)$ is not necessarily isomorphic to $X$ then understanding the autoequivalences of $\der$ alone does not give us an immediate relationship to the Picard group of $\boxtimes$. \\
We should think of the following corollary as a monoidal version of the Bondal-Orlov reconstruction theorem:
\begin{cor}\label{cor:monoidalbondalorlov}
Let $X$ be as above, then if $\omega_{X}[n]$ is an invertible object for a tensor triangulated structure $\boxtimes$ on $\der$ with unit $\Ox_{X}$ then $\boxtimes$ and $\dtee$ coincide on objects.
\begin{proof}
As $\omega_{X}$ is $\boxtimes$-invertible, Corollary \ref{cor:amplepicard} tells us that for any $A\in \der$ we have
\[ \omega_{X}\dtee A \cong \omega_{X}\boxtimes A. \]
But we can resolve any other complex $B$ in terms of derived powers of the canonical sheaf, by the exactness of $\boxtimes$ we have then 
\[ B\dtee A \cong B\boxtimes A. \]
\end{proof}
\end{cor}
The nature of this result comes precisely from the fact that the tensor triangulated category structure $(\boxtimes, \mathbbm{1})$ does not necessarily come from a derived equivalence $\der\simeq D^{b}(Y)$, while the extra assumption on the unit is required the result is in this direction slightly more general than the original theorem. \\
\begin{cor}\label{cor:isopicard}
Let $X$ be a variety with ample (anti-)canonical bundle, suppose $(\Ox_{X},\boxtimes)$ is a tensor triangulated structure on $\der[X]$ such that $Pic(\boxtimes)\cong Pic(\der[X])$ via the assignment $U\mapsto U$ then $\boxtimes$ coincides with $\dtee$ on objects.
\begin{proof}
In this case if this morphism is an isomorphism, then $\omega_{X}$ is $\boxtimes$-invertible and the result follows from the previous corollary.
\end{proof}
\end{cor}
In fact if we are under the same hypothesis for $X$ then as soon as we are able to show that the generators of $Pic(\der,\dteee)$ are $\boxtimes$-invertible then by the previous corollary there must be an equivalence between $\boxtimes$ and $\dtee$.
\begin{exmp}
Let $X=\mathbb{P}^{n}$ be the projective space, in this case we know that $Pic(D^{b}(X))=\mathbb{Z}\oplus \mathbb{Z}$ corresponding to the line bundles plus their shifts. The result above then says that whenever there is a tensor triangulated structure $\boxtimes$ on $D^{b}(X)$ with unit $\Ox_{X}$ then the Picard group of this tensor structure must necessarily be a subgroup of $\mathbb{Z}\oplus \mathbb{Z}$. \\
If $\omega_{X}=\Ox_{X}(-n-1)$ is $\boxtimes$-invertible then we get that $\boxtimes$ coincides with $\dtee$. Similarly if $\Ox_{X}(-1)$ is $\boxtimes$-invertible.
\end{exmp}
One natural question to ask when working with Picard groups of tensor triangulated category structures is what is the relationship with line bundles on the associated space. In \cite[Proposition 4.4]{balmer2007gluing} Balmer and Favi prove the following result:
\begin{prop}\label{prop:balmergluing}
Let $X$ be a scheme and consider $Perf(X)$ its derived category of perfect complexes. Then there is a split short exact sequence of abelian groups 
\[ 0\to Pic(X)\to Pic(Perf(X),\dtee)\to C(X;\mathbb{Z})\to 0\]
where $C(X;\mathbb{Z})$ stands for the group of locally constant functions from $X$ to $\mathbb{Z}$.
\end{prop}
Again under the hypothesis of $X$ having an ample (anti-)canonical bundle, by using Proposition \ref{prop:balmergluing} we see that for a TTC $(\boxtimes,\Ox_{X})$, if $Spc(\boxtimes)$ is a scheme then the Picard group of $Spc(\boxtimes)$ must be a subgroup of the Picard group of $Pic(\der, \boxtimes)\leq Pic(\der,\dtee)$. So a line bundle in $Spc(\boxtimes)$ has to be $\dtee$-invertible. \\ 
\begin{obs}\label{obs:reconstructiontensorantieu}
From Bondal-Orlov's reconstruction original proof we know that it is actually possible to fully characterize line bundles up to a shift from purely categorical properties. Given the importance of the Picard group of the variety, we can ask whether it is possible to reconstruct the derived tensor product in $\der$ without having to pass through a reconstruction theorem. \\
In \cite{143661} Antieau sketches a construction in which by considering invertible objects ( in the sense of Bondal and Orlov ) one can define a collection of tensor products $\dteee_{U}$ by exploiting the resolution by derived tensor powers of $\omega_{X}$. \\
The idea is to pick an invertible object $U$ which are shown in \cite{bondal2001reconstruction} to be isomorphic to a shift of a line bundle in $X$, then by use of the resolution we only need to define the products $\omega_{X}^{\dtee i}[ni]\dteee_{U}A^{\bullet}$ for any object $A^{\bullet}$. As the Serre functor $S\simeq \_\dtee \omega_{X}[n]$ comes with the categorical structure alone then we can set these products to be simply $S^{i}(A^{\bullet})$. \\
These tensor products $\dteee_{U}$ have $U$ as unit and all of them have $X$ as Balmer spectrum.
\end{obs}
In general for a triangulated category $\T$ we have an action by $Aut(\T)$ on the collection $TTS(\T)$.\\
If $(\otimes,\mathbbm{1})\in TTS(\T)$ and $\phi\in Aut(\T)$ we have a tensor structure defined by \[ X\otimes_{\phi} Y:= \phi^{-1}(\phi(X)\otimes \phi(Y)). \]
And with unit given by $\phi^{-1}(\mathbbm{1})$. \\
We have now justified enough the following definition:
\begin{defn}\label{defn:tens}
Let $\T$ be a triangulated category, denote by $TTS(\T)$ the collection of equivalence classes of tensor triangulated category structures on $\T$. Where we consider two tensor triangulated category structures to be equivalent if there is a monoidal equivalence between the two of them. \\
\end{defn}
To keep some control and avoid counting structures coming from autoequivalences as we saw, we should at least fix the unit object.
\begin{defn}
Let $\T$ be a triangulated category and $U\in \T$ an object. Then the set $TTS_{U}(\T)$ is the set of equivalence classes of tensor triangulated structures on $\T$ where $U$ is the unit.
\end{defn}
It is this set the one we are mainly interested in classifying. \\
Let us finish by discussing the original Bondal-Orlov reconstruction theorem in terms of the results we have shown so far. \\
\begin{thm}\label{thm:newbondalorlov}
Let $X$ be a variety with ample (anti-)canonical divisor, and let $\boxtimes$ be a tensor triangulated structure on $\der$ with unit $\Ox_{X}$. Suppose $Spc(\boxtimes)$ is a smooth projective space with ample (anti-)canonical bundle and that there is an equivalence $\der\simeq D^{b}(Spc(\boxtimes))$, then $X\cong Spc(\boxtimes)$
\begin{proof}
In fact the only thing to note here is that as $Spc(\boxtimes)$ has ample (anti-)canonical bundle then $\omega_{X}$ has to be $\boxtimes$-invertible. Indeed, we recall that one can pick the equivalence $\der\simeq D^{b}(Spc(\boxtimes))$ to send $\omega_{X}$ to $\omega_{Spc(\boxtimes})$ and then the assertion follows by applying Corollary \ref{cor:amplepicard} to $Spc(\boxtimes)$ we obtain that $Pic(\der, \dtee)$ has to be isomorphic via the assignment $\Ll\mapsto \Ll$  to a subgroup of $Pic(\der, \boxtimes)$. Since $\omega_{X}$ is $\boxtimes$-invertible, by Corollary \ref{cor:monoidalbondalorlov} we obtain our result. 
\end{proof}
\end{thm}
\begin{obs}
We need to explain our choice of hypothesis here. On the first hand the assumption that $Spc(\boxtimes)$ is a smooth projective variety is necessary just as in the original Bondal-Orlov theorem formulation. We have added a couple more assumptions, however. We suppose that the (anti-)canonical bundle of $Spc(\boxtimes)$ is also ample to highlight the use of the monoidal structures in the theorem. This hypothesis is however not necessary as it can be directly deduced from the derived equivalence between the two spaces, just as in the original proof of Bondal and Orlov. Alternatively, we can formulate the theorem as follows:
\begin{thm}
Let $X$ be a variety with ample (anti-)canonical divisor, and let $\boxtimes$ be a tensor triangulated structure on $\der$ with unit $\Ox_{X}$. Suppose $Spc(\boxtimes)$ is a smooth projective space, and that we have an equivalence $\der\simeq D^{b}(Spc(\boxtimes))$, then $X\cong Spc(\boxtimes)$.
\end{thm}
Of more importance is perhaps the choice of unit, as we have seen that there are tensor triangulated category structures on the derived category of such a variety which will produce very different spaces under the Balmer reconstruction. This choice of unit allows us to keep some control in the classification of structures producing the same space. \\
A natural next step for future work would be to deal with the possible sort of objects which can be units for such a structure. \\
\end{obs}
\begin{obs}
We wish to point out that there is some nuance in the way in which Bondal-Orlov follows from our results as we make use of some important technical results from the original proof. We expect however that the discussion in this work has provided enough of a justification and motivation for looking at this problem in terms of monoidal structures. 
\end{obs}
We can close our discussion with the following theorem:
\begin{thm}\label{thm:uniquetens}
Let $X$ be a smooth projective variety with (anti-)canonical bundle. Consider a tensor triangulated category structure $\boxtimes$ on $\der$ such that $\Ox_{X}$ is its unit and $Spc(\boxtimes)$ is isomorphic to $X$, then $\boxtimes$ and $\dtee$ coincide on objects.
\end{thm}
This however does not fully classify $TTS_{\Ox_{X}}(\der)$ as we require Balmer's spectrum to be a Fourier-Mukai partner, but there is no reason to expect in general a relationship between the derived category of the Balmer spectrum and the original triangulated category. \\
The lack of morphisms between a space $X$ and the Balmer spectrum $Spc(\boxtimes)$ for some tensor triangulated structure, and thus of functors between the derived categories of these two spaces is one of the obstacles to being able to understand the possible structures $\boxtimes$.
\hypersetup{linkcolor=red}
\bibliographystyle{alpha}
\bibliography{Main}

\end{document}